\documentclass[12pt]{article}
\usepackage{amsmath}
\usepackage{bbm}
\usepackage{mathbbold}
\usepackage{amsfonts}
\usepackage{cite}
\usepackage{mathrsfs}
\usepackage{amsfonts,amssymb}

\numberwithin{figure}{section}
 \numberwithin{equation}{section}

\newtheorem{theorem}{Theorem}[section]
\newtheorem{proposition}[theorem]{Proposition}
\newtheorem{definition}[theorem]{Definition}
\newtheorem{corollary}[theorem]{Corollary}
\newtheorem{lemma}[theorem]{Lemma}
\newtheorem{remark}[theorem]{Remark}

\renewcommand{\theequation}{\thesection.\arabic{equation}}

\newcommand{\cH}{{\mathcal H}}

\newcommand{\cT}{{\mathcal T}}

\newcommand{\cM}{{\mathcal M}}

\newcommand{\sA}{{\mathscr A}}

\newcommand{\sX}{{\mathscr X}}

\def\R{\mathbb{R}}

\def\be{\begin{equation}}
\def\ee{\end{equation}}
\def\ba{\begin{array}}
\def\ea{\end{array}}
\def\benu{\begin{enumerate}}
\def\eenu{\end{enumerate}}
\def\bt{\begin{theorem}}
\def\et{\end{theorem}}
\def\bl{\begin{lemma}}
\def\el{\end{lemma}}
\def\br{\begin{remark}}
\def\er{\end{remark}}
\def\bd{\begin{definition}}
\def\ed{\end{definition}}
\def\bp{\begin{proposition}}
\def\ep{\end{proposition}}
\def\bc{\begin{corollary}}
\def\ec{\end{corollary}}

\def\b{\beta}
 \def\de{\delta}  \def\nab{\nabla}
\def\lam{\lambda} 
 \def\sig{\sigma}
\def\w{\omega}\def\W{\Omega}
\def\a{{\alpha}}
\def\va{\varphi}

\def\0{\theta}

\def\<{\subset}

\def\.{\cdot}

\def\A{\forall}
\def\ol{\overline}

\def\Cap{\bigcap}
\def\Cup{\bigcup}

\def\ra{\rightarrow}
\def\stac{\stackrel}
\def\~{\stac{\sim}}
\def\8{\infty}
\def\X{\times}

\def\mb{\mbox}

\def\p{\rho}
\def\Hs{\hspace{0.8cm}}
\def\hs{\hspace{0.4cm}}
\def\Vs{\vskip10pt}
\def\vs{\vskip5pt}

\def\mB{\mbox{B}}

\def\[{\left[}
\def\]{\right]}
\def\({\left(}
\def\){\right)}
\def\square{\qed}

\begin{document}
\begin{center}
{\bf\Large {An Invariant Set Bifurcation Theory for Nonautonomous Nonlinear Evolution Equations \footnote{This work was supported by the National Natural Science Foundation of China [11871368]. }}}
\end{center}

\begin{center}
Xuewei Ju,\Hs Ailing Qi
\end{center}

\begin{center}
{\footnotesize { Department of Mathematics, Civil Aviation University of China\\
     Tianjin 300300,  China }}

\end{center}
\footnote{{\em E-mail}: xwju@cauc.edu.cn (X.W. Ju), alqi@cauc.edu.cn (A.L. Qi).}

{\bf Abstract.}
In this paper we establish an invariant set bifurcation theory  for the nonautonomous dynamical system $(\va_\lam,\0)_{X,\cH}$ generated by the evolution equation \be\label{e0}u_t+Au=\lam u+p(t,u),\hs p\in \cH=\cH[f(\.,u)]\ee on a Hilbert space $X$, where $A$ is a sectorial operator, $\lam$ is the bifurcation parameter, $f(\.,u):\R\ra X$ is translation compact, $f(t,0)\equiv0$ and $\cH[f]$ is the hull of $f(\.,u)$. Denote by $\va_\lam:=\va_\lam(t,p)u$ the cocycle semiflow generated by the equation.
Under some other assumptions on $f$, we show that as the parameter $\lam$ crosses an eigenvalue $\lam_0\in\R$ of $A$, the system bifurcates from
 $0$  to a nonautonomous invariant set $B_\lam(\.)$ on one-sided neighborhood of $\lam_0$. Moreover, $$\lim_{\lam\ra\lam_0}H_{X^\a}\(B_\lam(p),0\)=0,\hs p\in P,$$
where  $H_{X^\a}(\.,\.)$ denotes  the Hausdorff semidistance in $X^\a$  (here $X^\alpha$ ($\a\geq0$) defined below is the  fractional power spaces associated with $A$).

 Our result is based on the pullback attractor bifurcation on the local central invariant manifolds $\cM^\lam_{loc}(\.)$.
 \Vs{\bf Keywords.} Stability of pullback attractors; local invariant manifolds; nonautonomous invariant set bifurcations.
\section{Introduction}
Invariant set bifurcation theory of autonomous dynamical systems has been  extremely well developed \cite{Kie,Kra,MW1,MW2,MW3,MW4,AY,Ward,chow1,Nirenberg,San,Rybakowski,LW}.
A relatively simpler but important case is that of bifurcations from equilibria, including bifurcation to multiple equilibria (static bifurcation) and to periodic solutions (Hopf bifurcation)  (see
among others, \cite{chow1,Nirenberg}).  Ma and Wang \cite{MW1} and Sanjurjio \cite{San} developed a local attractor bifurcation theory.
Roughly speaking, if the trivial equilibrium $e$ of an autonomous system changes from an attractor to a repeller on the local center manifold of the equilibrium when the bifurcation parameter $\lam$ crosses a critical value $\lam_0$, then the system bifurcates a compact invariant set   $K$ which is an attractor of the system restricted to the center manifold.
Chow and Hale \cite{chow1}  started to  discuss stability and bifurcation phenomena
associated with more general invariant sets, e.g. periodic orbits. Using Conley index theory, Rybakowski \cite{Rybakowski} and  Li and  Wang \cite{LW} developed  global bifurcation theorems to discuss bifurcation phenomena of nonlinear autonomous evolution equations.

However, except some relatively simple nonautonomous cases, there are  few papers studying the invariant set bifurcation for nonautonomous dynamical system. In \cite{LRS} Langa et al. presented a collection of examples to illustrate  bifurcation phenomena in nonautonomous ordinary
differential equations. In \cite{CLR}  Carvalho et al. studied the structure of the pullback attractor
for a nonautonomous version of the Chafee-Infante equation, and
investigated the bifurcations that this attractor undergoes as  bifurcation parameter varies.

Unlike autonomous dynamical systems for which forward dynamics is studied, pullback dynamics is much more
natural than the more familiar forward dynamics for nonautonomous dynamical systems.  This makes it very difficult to extend the invariant set bifurcation theory of autonomous systems to nonautonomous systems. Our approach in the paper is to treat the nonautonomous
system as a cocycle semiflow over a suitable base space rather than a process. The biggest advantage of the cocycle semiflow framework is that in many cases the base spaces are compact, while the default base space $\R$ (real number set) for processes is unbounded. Based on the compactness of the base spaces, we can establish  the equivalence between pullback attraction of cocycle semiflow and forward attraction of the associated autonomous semiflow. This device makes the dynamics of such a nonautonomous system appear like those of an autonomous system.

Without the compactness assumption on the base spaces, the upper semicontinuity of global pullback attractors for nonautonomous systems was obtained in Caraballo and Langa \cite{CL}. However,  compact forward invariant sets of the perturbed systems are required to guarantee the existence of  perturbed pullback attractors. In the paper, we suppose that the base spaces of cocycle semiflows considered are compact, which will require some
restrictions on the nonlinearities. As a result, after  introducing the notion of  (local) pullback attractors (see Definition \ref{de1}), we can establish a general result on the stability of local pullback attractors as the perturbation parameter is varied. Based on this result, a local pullback attractor bifurcation theory can be developed.  This  can be regarded as a nonautonomous  generalization of autonomous attractor bifurcation theory in \cite{MW1}. Finally, we study the bifurcation of invariant sets for the cocycle semiflow $\va_\lam$ generated by the nonautonomous nonlinear evolution equation \eqref{e0}. We first construct a local central invariant manifold $\cM^\lam_{loc}(\.)$ for  $\va_\lam$ with $\lam$ near $\lam_0$. Under further assumptions on $f$ to ensure that $0$ is a pullback attractor for $\va_{\lam_0}$, we then restrict $\va_\lam$ to  $\cM^\lam_{loc}(\.)$ and obtain a pullback attractor bifurcation on  $\cM^\lam_{loc}(\.)$ as $\lam$ crosses $\lam_0$. It leads to an invariant set bifurcation for $\va_\lam$. It is worth mentioning that if $0$ is not an attractor but a repeller for  $\va_{\lam_0}$, our result still holds. Denote by $B_\lam(\.)$ the bifurcated invariant set. We further know that $$\lim_{\lam\ra\lam_0}H_{X^\a}\(B_\lam(p),0\)=0,\hs p\in P.$$

This paper is organized as follows. In Section 2, we present respectively some basic facts in autonomous and nonautonomous dynamical systems which will be required in the rest of the work. Section 3 deals with the
 stability of pullback attractors as bifurcation parameter varies. In Section 4, we establish an invariant set bifurcation theory for \eqref{e0}. We illustrate the main results with an example in Section 5. Finally, Section 6 contains the proofs of two propositions presented earlier in the paper.

 \section{Preliminaries}
In this section we introduce  some basic definitions and notions \cite{CKS, CV}.

Let $X$ be a complete metric space with metric $d(\.,\.)$.
Given $M\subset X$,
we denote $\ol M$ and int$\,M$ the closure and interior  of any subset  $M$ of $X$,  respectively.
A set $U\subset X$ is called a
neighborhood of $M\subset X$, if $\overline{M}\subset \mbox{int}\,U$. For any $\p>0$, denote by $$\mb{B}_{X}(M,\p):=\{x\in X:\,\,d(x,M)<\p\}$$ the $\p$-neighborhood of $M$ in $X$, where $d(x,M)=\inf_{y\in M}d(x,y)$.

The  Hausdorff semidistance  in $X$ is defined  as $$H_{X}(M,N)=\sup_{x\in M}d(x,N),\Hs\A\,M,N\subset X.$$

\subsection{Semiflows and attractors}Let $\R^+=[0,\8)$.
A continuous mapping $S:\R^+\times X\ra X$ is called a semiflow on $X$, if it satisfies
\begin{enumerate}
\item[$i)$]$S(0,x)=x$ for all $x\in X$; and
\item[$ii)$]$S(t+s,x)=S(t,S(s,x))$ for all $x\in X$ and  $t,s\in\R^+$.
\end{enumerate}

\vs

Let $S$ be a given semiflow on $X$. As usual, we will rewrite $S(t,x)$ as $S(t)x$.

A set $B\<X$ is called {\bf invariant} (resp. {\bf positively invariant}) under $S$ if $S(t)B=B$ (resp. $S(t)B\<B$) for all $t\geq0$.

Let $B$ and $C$ be subsets of $X$. We say that $B$ {\bf attracts} $C$ under $S$, if $$\lim_{t\ra\8}H_X(S(t)C,B)=0.$$

\bd
A compact subset $\sA\<X$ is called an {\bf attractor} for $S$, if it is invariant under $S$ and attracts
one of neighborhood of itself.
\ed
It is well known that if $U$ is a compact positively invariant set of $S$, then the omega-limit
set $\w(U):=\bigcap_{T\geq0}\overline{\bigcup_{t\geq
T}S(t)U}$ is an attractor of $S$. The definition of the attraction basin of the attractor and other properties of local attractors can be found in \cite{Con,MM,Rybakowski}.
%




\subsection{Cocycle semiflows and pullback attractors}
A nonautonomous system consists of a ``base flow'' and a ``cocycle semiflow'' that is in some sense driven by the base flow.

A {\bf base flow} $\{\0_t\}_{t\in\R}:=\{\0(t)\}_{t\in\R}$  is a flow on a metric space $P$ such that
$\0_tP=P$ for all $t\in\R$.

\bd
A {\bf cocycle semiflow} $\va$ on the phase space $X$ over $\0$ is a continuous mapping $\va:\R^+\times P\times X\ra X$ satisfying \begin{enumerate}
\item[$\bullet$] $\va(0,p,x)=x$,
\item[$\bullet$]$\va(t+s,p,x)=\va(t,\0_sp,\va(s,p,x))$ (cocycle
property). \end{enumerate}
\ed
\br If we replace $\R^+$ by $\R$ in the above definition, then $\va$ is called a cocycle flow on $X$.\er

We usually denote $\va(t,p)x:=\va(t,p,x)$. Then $\{\va(t,p)\}_{t\geq0,\,p\in P}$ can be viewed as a family of continuous mappings on $X$.

For convenience in statement, a family of subsets $B(\.):=\{B_p\}_{p\in P}$ of $X$ is called a {\bf nonautonomous set} in $X$.
As usual, we will rewrite  $B_p$ as $B(p)$, called the $p$-{\bf section} of $B(\.)$. We also denote $\mathbb{B}$ the union of the sets $B(p)\X\{p\}$ ($p\in P$), i.e., $$\mathbb{B}=\Cup_{p\in P}B(p)\X\{p\}.$$
Note that  $\mathbb{B}$ is a subset of $X\X P$.

A nonautonomous set $B(\.)$ is said to be closed (resp. open, compact), if each section $B(p)$ is closed (resp. open, compact) in $X$.
A nonautonomous set $U(\.)$ is called a {\bf
neighborhood} of $B(\.)$, if $\ol B(p)\subset \mbox{int}\,U(p)$ for each $p\in P$.

A nonautonomous set $B(\.)$ is said to be {\bf invariant} (resp. {\bf forward invariant}) under $\va$ if for $t\geq0$, $$\va(t,p)B(p)=B(\0_tp),\Hs p\in P.$$
$$(\hbox{resp. } \va(t,p)B(p)\<B(\0_tp),\Hs p\in P.)$$
Let $B(\.)$ and $C(\.)$ be two nonautonomous subsets of $X$.  We say that $B(\.)$ {\bf pullback  attracts} $C(\.)$ under $\va$ if for any $p\in P$,$$\lim_{t\ra\8}H_{X}(\va(t,\0_{-t}p)C(\0_{-t}p),\,B(p))=0.$$

Let $\va$ be   a given cocycle semiflow  on $X$ with driving system $\0$ on base space  $P$. The (autonomous) semiflow  $\Phi:=\{\Phi(t)\}_{t\geq0}$ on $Y:=P\X X$, given by $$\Phi(t)(p,x)=(\0_t p,\phi(t,p)x),\hs t\geq0,$$ is called the {\bf  skew product semiflow} associated to $\va$.
The following fundamental result studies the relationship between the pullback attraction of $\va$ and attraction of $\Phi$.  The proof is given in Appendixes.
\bp\label{P1}Let $(\va,\0)_{X,P}$ be a nonautonomous system, and  let $\Phi$ be the skew-product flow associated to $\va$. Let $K(\.)$ and $B(\.)$ be two nonautonomous sets.  Suppose $P$ and $K_P:=\overline{\bigcup_{p\in P}K(p)}\<X$ are both compact. Then $K(\.)$ pullback attracts $B(\.)$ through $\va$ if and only if $\mathbb{K}:=\bigcup_{p\in P}K(p)\X \{p\}$ attracts $\mathbb{B}:=\bigcup_{p\in P}B(p)\X \{p\}$ through $\Phi$.\ep

\bd\label{de1}
Let $(\va,\0)_{X,P}$ be a nonautonomous system. A nonautonomous set $A(\.)$ is called a {\bf (local) pullback attractor} for $\va$ if it is compact, invariant and  pullback attracts a neighborhood $U(\.)$ of itself.
\ed
The local pullback attractor defined here, very similar to the notion of a past attractor in  Rasmussen \cite{Ras}, can be seen as a nature nonautonomous generalization of the local attractor from the
autonomous theory.
Similar to the case of autonomous systems, if $U(\.)$ is a compact forward invariant set of $\va$, then  the omega-limit
set $\w(U)(\.)$ defined as
$$\w(U)(\w)=\bigcap_{T\geq0}\overline{\bigcup_{t\geq
T}\varphi(t,\0_{-t}\w)U(\0_{-t}\w)},\hs \w\in \W$$ is a pullback attractor of $\va$.
For instance, consider the following simple system on $X=\R$:
\be\label{eq2.1}x'(t)=-3x+p(t)x^3,\hs p\in\cH[h],\ee
where $h(t)=2+\sin t$ and $\cH[h]$ is its hull which is the closure for the  uniform convergence topology of the set of $t$-translates of $h$. The translation map $\0_t:\cH\ra\cH$ given by $\0_tp(s)=p(t+s)$ defines a flow on $\cH$. Then the unique solution of \eqref{eq2.1} define a cocycle flow on $X$ given by $\va(t,p)x_0=x(t,0;p;x_0)$. Since $$\frac{1}{2}\frac{d}{dt}x^2=-3x^2+p(t)x^4\leq -3(x^2-x^4)<0$$ provided that $|x|\leq1/2$. Therefore $[-1/2,1/2]$ is a forward invariant set of $\va$ and it is pullback attracted by the pullback attractor $0$.
It is worth noting that $0$ is only a local pullback attractor. Indeed, $$\frac{1}{2}\frac{d}{dt}x^2=-3x^2+p(t)x^4\geq -3x^2+x^4>0$$ provided that $|x|\geq2$. It follows that $0$ is only a local pullback attractor of $\va$.

In general, it is difficult to define the attraction basin of a pullback attractor. Fortunately, under the assumptions of Proposition \ref{P1}, we can define the pullback attraction basin of a pullback attractor $A(\.)$. Specifically, we have
\bd Let $(\va,\0)_{X,P}$ be a nonautonomous system, and  let $\Phi$ be the skew-product flow associated to $\va$. Suppose $P$ is compact. Let $A(\.)$ be a pullback attractor of $\va$ such that $A_P:=\overline{\bigcup_{p\in P}A(p)}$ is compact. Let $B(\mathbb{A})=\{(x,p):\mathbb{A} \hbox{ attracts } (x,p) \hbox{ through } \Phi \}$ be the attractor basin of $\mathbb{A}$ under $\Phi$. Then the {\bf pullback  attraction basin} $B(A)(\.)$ of $A(\.)$ can be defined as $$B(\mathbb{A})=\bigcup_{p\in P}B(A)(p)\X\{p\}.$$\ed

\section{Stability of pullback attractors}
We now establish a result on the stability of pullback attractors under a small perturbation. In fact,
we prove a continuity result with respect to the Hausdorff semidistance.

Let $X$ be a Banach space with norm $\|\.\|$, and let $A$ be a sectorial operator on $X$. Pick a number $a>0$ such that $$\hbox{Re}\,\sig(A+aI)>0.$$
Denote  $\Lambda=A+aI$. For each $\a\geq0$, define the
fractional power space as $X^\a=D(\Lambda^\a)$, which is equipped with the norm $\|\.\|_\a$ defined by
$$\|x\|_\a=\|\Lambda^\a x\|,\hs x\in X^\a.$$
Note  that the definition of $X^\a$ is independent of the choice of the number  $a$.
If $A$ has compact resolvent, the inclusion $X^{\a'}\hookrightarrow X^\a$ is compact for $\a'>\a\geq0$.

Let $\va_{\lam_0}$ ($\lam_0\in\R$) be   a given cocycle semiflow  on  $X$ with driving system $\0$ on base space  $P$. For $\de>0$, denote $I_{\lam_0}(\de):=(\lam_0-\de,\lam_0+\de)$. Assume that $\va_\lam$, $\lam\in I_{\lam_0}(\de)$ is a small
perturbation of the given flow $\va_{\lam_0}$ based on $P$.
Let us make the following assumptions:
\benu
 \item[\,\,{\bf (H1)}:] The base space $P$ is compact.
 \eenu
 \benu  \item[\,\,{\bf (H2)}:]
For every  $T>0$ and compact subset $B$ of $X$, we have
\be\label{e8}\lim_{\lam\ra\lam_0}\|\va_\lam(t,p)x-\va_{\lam_0}(t,p)x\|_\a=0,\ee
uniformly with respect to $(t,x)\in [0,T]\X B$ and $p\in P$.
\eenu

Under the assumptions {\bf (H1),\,(H2)}, we can get a result on the stability of pullback  attractors.
\bt\label{th1} Let $A_{\lam_0}(\.):=\{A_{\lam_0}(p)\}_{p\in P}$ be an attractor of the cocycle semiflow $\va_{\lam_0}$  which pullback attracts a neighborhood $U(\.)$ of itself. Let  $$\mathbb{U}:=\bigcup_{p\in P}U(p)\X \{p\}\hs\hbox{and} \hs \mathbb{A}_{\lam_0}:=\overline{\bigcup_{p\in P}A_{\lam_0}(p)\X \{p\}}.$$ Assume $\mathbb{U}$ is a compact neighborhood of $\mathbb{A}_{\lam_0}$ in $Y=X\X P$, then under the assumptions {\bf (H1),\,(H2)}, the following statements hold.\begin{enumerate}
\item[$(\mathbbm{a})$]  There exists a small $\de>0$ such that for each $\lam\in I_{\lam_0}(\de)$, $\va_\lam$ has a pullback attractor $A_\lam(\.)$ such that \be\label{eq1}\lim_{\lam\ra\lam_0}H_X(A_\lam(p),\bigcup_{p\in \cH}A_{\lam_0}(p))=0.\ee
    \item[$(\mathbbm{b})$] In addition, if $U(\.)$ is forward invariant, then \be\label{eq2}\lim_{\lam\ra\lam_0}H_X(A_\lam(p),A_{\lam_0}(p))=0.\ee\end{enumerate}\et
{\bf Proof.} $(\mathbbm{a})$  By the compactness of $\mathbb{U}$, we know that ${A_{\lam_0}}_P:=\overline{\bigcup_{p\in P}A_{\lam_0}(p)}$ is compact. Since $A_{\lam_0}(\.)$ pullback attracts $U(\.)$ and $P$ is compact, by Proposition  \ref{P1},
$\mathbb{A}_{\lam_0}$ attracts $\mathbb{U}$ through $\Phi_{\lam_0}$. Since $\mathbb{U}$ is a neighborhood of $\mathbb{A}_{\lam_0}$, one knows that $\mathbb{A}_{\lam_0}$ is an attractor of $\Phi_{\lam_0}$.  By the assumption {\bf (H2)}, for any compact set $B\subset X$, we have that \be\lim_{\lam\ra\lam_0}H_Y(\Phi_\lam(t)(x,p),\Phi_{\lam_0}(t)(x,p))=\lim_{\lam\ra\lam_0}\|\va_\lam(t,p)x-\va_{\lam_0}(t,p)x\|_\a=0\ee uniformly with respect to $t\in[0,T]$ and $(x,p)\in B\X P$.
Then by the stability of the autonomous attractors \cite[Theorem 4.1]{LZ},  there exists a $\de>0$ (independent of $p\in P$) such that for each $\lam\in I_{\lam_0}(\de):=(\lam_0-\de,\lam_0+\de)$, $\Phi_{\lam}$
 has an attractor $\mathbb{A}_\lam$ contained in $\mathbb{U}$. Moreover,  \be\label{eeq3}\lim_{\lam\ra\lam_0}H_Y(\mathbb{A}_\lam,\mathbb{A}_{\lam_0})=0.\ee
Write $\mathbb{A}_\lam$ as $\bigcup_{p\in P}A_\lam(p)\X \{p\}$, $\lam\in I_{\lam_0}(\de)$.
Using Proposition \ref{P1} again, we have that $A_\lam(\.)$  pullback attracts $U(\.)$ through $\va_\lam$, i.e., $A_\lam(\.)$ is a pullback attractor of $\va_\lam$. \eqref{eq1} is a direct consequence of \eqref{eeq3}.

To complete the proof of $(\mathbbm{b})$, we shall prove \eqref{eq2} by contradiction. Thus, let us assume
that there exist $\sigma>0$ and  sequences $\lam_j\ra \lam_0$, as $j\ra\8$, $x_j\in A_{\lam_j}(p)$ such that
\be\label{eq4}d_X(x_j,x)>\sig,\hs\hbox{ for all }x\in A_{\lam_0}(p).\ee Note that $$x_j=\va_{\lam_j}(n,\0_{-n}p)x_j^n,\hs \hbox{ for some }x_j^n\in A_{\lam_j}(\0_{-n}p).$$
Similar to the argument in $(\mathbbm{a})$, we can assume that $A_{\lam_j}(p)\<U(p)$, thus $x_j\in U(p)$. By the compactness of $U(p)$, there exists a subsequence of $x_j$ (still denoted by $x_j$) which converges to some $x_0\in U(p)$. Now, for each fixed $n$ we have $x_j^n\in U(\0_{-n}p)$ so that
there is a further subsequence of $x_j^n$ (still denoted by $x_j^n$) which converges to some $x_0^n\in U(\0_{-n}p)$.
On the other hand, for any given $\nu>0$, we can use the assumption {\bf (H2)} and the continuity of $\va(n,\0_{-n}p)$ to show that for $j$ large enough,
\begin{equation*}
\begin{split}&d(\va_{\lam_j}(n,\0_{-t}p)x_j^n, \va_{\lam_0}(n,\0_{-t}p)x_0^n)\\
\leq&d(\va_{\lam_j}(n,\0_{-t}p)x_j^n, \va_{\lam_0}(n,\0_{-t}p)x_j^n)\\
&+d(\va_{\lam_0}(n,\0_{-t}p)x_j^n, \va_{\lam_0}(n,\0_{-t}p)x_0^n)\\
\leq&\nu+\nu.\end{split}
\end{equation*}
Then, for each fixed $n\in \mathbb{N}$,
$$x_0=\lim_{j\ra \8}x_j=\lim_{j\ra \8}\va_{\lam_j}(n,\0_{-n}p)x^n_j=\va_{\lam_0}(n,\0_{-n}p)x_0^n.$$
Since $U(p)$ is forward invariant, we have $$x_0\in \bigcap_{n\in \mathbb{N}}\va_{\lam_0}(n,\0_{-n}p)U(\0_{-n}p)=A_{\lam_0}(p),$$which contradicts \eqref{eq4}.
The proof is complete. $\Box$

The main contribution of Theorem \ref{th1} is the existence of pullback attractor $A_\lam(\.)$ for $\va_\lam$ as $\lam$ near $\lam_0$, while the argument of the upper semicontinuity of pullback attractors is an adaptation of that of \cite{CL}.

The conditions of the following results may be easier to be verified in applications.
\bc\label{c1} Let $A_{\lam_0}(\.):=\{A_{\lam_0}(p)\}_{p\in P}$ be an attractor of the cocycle semiflow $\va_{\lam_0}$ and $U\<X$ is a compact forward invariant neighborhood of $A_{\lam_0}(\.)$. Then under the assumptions {\bf (H1),\,(H2)}, there exists a small $\de>0$ such that for each $\lam\in I_{\lam_0}(\de)$, $\va_\lam$ has a pullback attractor $A_\lam(\.)$ satisfying $$\lim_{\lam\ra\lam_0}H_X(A_\lam(p),A_{\lam_0}(p))=0.$$ \ec

\section{Invariant set bifurcation for nonautonomous nonlinear evolution equations}
Based on the general result of the stability of pullback attractors, in the section we can establish some results on invariant set bifurcation for nonautonomous dynamical systems.

\subsection{Problem and mathematical setting}
From now on, we assume $X$ is a Hilbert space with inner product $(\.,\.)$. We will consider and study invariant set bifurcation of the evolution  equation
\be\label{eeq1}u_t+A u=\lam u +f(t,u)\ee on $X$,
where $\lam\in \R$ is a bifurcation parameter,  the nonlinearity
 $f:\R \X X^\a\ra X$  is bounded continuous mapping satisfying \benu\item[{\bf (F1)}]  \be\label{e3}f(t,u)=o(\|u\|_\a),\hs \hbox{ as }\|u\|_\a\ra0\ee uniformly on $t\in\R$.
Moreover, there is $\b>0$ such that \be\label{eq1.7}\((f(t,u),u\)\leq-\b\.\kappa(u)\ee
 for $t\in \R$ and  $u\in X^\a$, where $\kappa:X\ra \R^+$ is a nonnegative function satisfying that $\kappa(u)=0$ if and only if $u=0$.
 \eenu Denote $k(\rho)$ the Lipschitz constant of $f(t,\.)$ in $\overline{\mB}_{X^\a}(\rho)$. Then by \eqref{e3}, $$\lim_{\rho\ra0}k(\rho)=0$$
and \be\label{eqq6}\|f(t,u_1)-f(t,u_2)\|\leq k(\rho)\|u_1-u_2\|_\a,\Hs \A\, u_1,u_2\in \overline{\mB}_{X^\a}(\rho).\ee

Denote $C_b(\R,X)$ the set of bounded continuous functions from $\R$ to $X$.  Equip $C_b(\R,X)$ with either the uniform convergence topology generated by the metric $$r(h_1,h_2)=\sup_{t\in\R}\|h_1(t)-h_2(t)\|,$$ or the compact-open topology generated by the metric
$$ r(h_1,h_2)=\sum_{n=1}^\8\frac{1}{2^n}\.\frac{\max_{t\in[-n,n]}\|h_1(t)-h_2(t)\|}{1+\max_{t\in[-n,n]}\|h_1(t)-h_2(t)\|}.$$
Then $C_b(\R,X)$ is a complete metric space.

Let $f(\.,u)\in C_b(\R,X)$ be the function in  \eqref{eeq1}. Define the hull of $f(\.,u)$ as follows $$\cH:=\cH[f(\.,u)]=\overline{\{f(\tau+\cdot,u):\,\,\tau\in\R\}}\,^{C_b(\R,X)}.$$
In application, $f(\.,u)$ is often taken as a periodic function, quasiperiodic function, almost periodic function, local almost periodic function \cite{CV, LD} or uniformly almost automorphic function \cite{YS}. In this case, the hull $\cH$ is a {\bf compact} metric space.
Accordingly, the translation group $\theta$ on $\cH$ is given by  $$\0_t p(\.,u)=p(t+\.,u),\Hs t\in\R,\,\,p\in\cH.$$
Instead of \eqref{eq1},  we will  consider the more general cocycle system in $X^\a$ (where $\a\in[0,1)$): \be\label{e1}u_t+Au=\lam u+p(t,u),\Hs p\in\cH.\ee

\bp\label{p3}\cite{D.H}\,
Let $A$ and $p$ be given as above. Assume that   $p$ is locally H$\ddot{o}$lder continuous in $t$. Then for each  $u_0\in X^\a$, there is a $T>t_0$ such that $(\ref{e1})$ has a unique solution $u(t)=u_\lam(t,t_0;u_0,p)$ on $[t_0,T)$ satisfying
\be\label{eeq2}u(t)=e^{-A(t-t_0)}x_0+\int_{t_0}^te^{-A(t-s)}[\lam u(s)+p(s,u(s))]ds,\hs t\in[t_0,T).\ee
\ep

For convenience, from now on we always assume that the unique solution $(\ref{eeq2})$ is globally defined.
Define
$$\va_\lam(t,p)u:=u_\lam(t,0;u,p),\Hs u\in X^\a,\,\,p\in\cH.$$
Then $\va_\lam$ is  a {\bf cocycle semiflow}  on $X^\a$ driven by the base flow $\0$ on  $\cH$.
Note that for each $p\in\cH$, $u(t)$ is a $p$-solution of $\va_\lam$ on an interval $J$ {\em if and only if} it solves the equation \eqref{e1} on $J$.

\subsection{Local invariant manifolds}
Let $\lam_0\in\R$ be an isolated eigenvalue of $A$. Suppose that
\benu\item[{\bf (F2)}] there is a $\eta>0$ such that the spectrum $$\sig(A)\cap\{z\in\mathbb{C}:\lam_0-\eta<\hbox{Re}z<\lam_0+\eta\}=\lam_0.$$ \eenu

Denote $A_{\lam}:=A-\lam$. Then for $\lam\in I_{\lam_0}(\eta/4):=(\lam_0-\eta/4,\lam_0+\eta/4)$, the spectrum
   $\sig(A_\lam)$ has a decomposition $\sig(A_{\lam})=\sig_c\cup \sig_+\cup\sig_-$, where $$\sig_c=\{\lam_0-\lam\},\hs\sig_+=\sig(A_{\lam})\cap\{\hbox{Re}\, \lambda>0\}\hs\hbox{ and }\hs\sig_-=\sig(A_{\lam})\cap\{\hbox{Re}\, \lambda<0\}.$$
Accordingly, the space $X$ has a direct sum
decomposition: $X=X_c\oplus X_+\oplus X_-$.
Denote $X_\pm=X_+\oplus X_-$ and $$X^\a_i:=X_i\cap X^\a,\Hs  i=c,+,-,\pm.$$
Note that $X_c^\a$ is finite dimensional.

Under the assumptions on $A$ and $f$, we can construct a local  invariant manifold for $\va_\lam$, $\lam\in I_{\lam_0}(\eta/8)$.
\bp \label{p2} Suppose  the assumptions {\bf (F1),\,(F2)} hold. Then there exists $\varrho>0$ such that the cocycle semiflow $\va_\lam$, $\lam\in I_{\lam_0}(\eta/8)$
 has a local invariant manifold $\cM^\lam_{loc}(\.):=\{\cM^\lam_{loc}(p)\}_{p\in\cH}$ in $X^\a$ which is
represented as $$\cM^\lam_{loc}(p)=\{y+\xi^\lam_p(y):y\in \overline{\mb{B}}_{X^\a_c}(\varrho)\},$$ where $\xi^\._p(\.):I_{\lam_0}(\eta/8)\X \overline{\mb{B}}_{X^\a_c}(\varrho)\ra X^\a_\pm$ is a  Lipschitz continuous mapping satisfying that
\be\label{eq1.5}\xi_p^{\lam}(0)=0\hbox{ and }\|\xi_p^{\lam}(y)-\xi_p^{\lam}(z)\|_\a\leq L_1\|y-z\|_\a\ee
 and \be\label{eq1.5'}\|\xi_p^{\lam_1}(y)-\xi_p^{\lam_2}(y)\|_\a\leq L_2|\lam_1-\lam_2|,\ee where $L_1>0$ is independent of $p\in P$ and $\lam\in I_{\lam_0}(\eta/8)$, and $L_2>0$ is independent of $p\in P$ and $y\in \overline{\mb{B}}_{X^\a_c}(\varrho)$.
\ep
The long proof of the above proposition is given in Appendixes.
\subsection{Invariant set bifurcation}
Firstly, let us restrict the equation \eqref{e1} on the invariant manifold $\cM_{loc}^\lam(\.)$, $\lam\in I_{\lam_0}(\eta/8)$.
Specifically, we study the finite dimensional equation  \be\label{eq3}y_t+(\lam_0-\lam)y=p(t,y+\xi^\lam_{\0_tp}(y)),\hs y\in \overline{\mb{B}}_{X^\a_c}(\varrho),\,p\in\cH.\ee
Denote  $\phi_\lam$   the cocycle flow  on $\overline{\mb{B}}_{X^\a_c}(\varrho)$ with driving system $\0$ on the base space  $\cH$  generated by \eqref{eq3}.

We first say that the condition {\bf (H2)} (in Section 3) holds for the cocycle flow $\phi_\lam$, $\lam\in I_{\lam_0}(\eta/8)$.  Specifically, we have the following result.
\bl\label{le3} For every  $T>0$, we have
\be\label{ee8}\lim_{\lam\ra\lam_0}\|\phi_\lam(t,p)y-\phi_{\lam_0}(t,p)y\|_\a=0,\ee
uniformly with respect to $(t,y)\in [0,T]\X \overline{\mb{B}}_{X^\a_c}(\varrho)$ and $p\in P$.\el
{\bf Proof.} For $\lam\in I_{\lam_0}(\eta/8)$, denote $y_\lam(t):=\phi_\lam(t,p)y$ and $v(t)=y_\lam(t)-y_{\lam_0}(t)$, then $v$ satisfies
\be\label{eq8}v_t+(\lam_0-\lam) y_\lam=p(t,y_\lam+\xi^{\lam}_{\0_tp}(y_\lam))-p(t,y_{\lam_0}+\xi^{\lam_0}_{\0_tp}(y_{\lam_0})).\ee
Note that $ \|y_\lam\|\leq\p$ and \begin{equation}\label{eq1.4}\begin{split}&\|p(t,y_\lam+\xi^{\lam}_{\0_tp}(y_\lam))-p(t,y_{\lam_0}+\xi^{\lam_0}_{\0_tp}(y_{\lam_0}))\|\\
\leq&k(\p)\big((L_1+1)\|v\|_\a+L_2|\lam-\lam_0|\big)\\
\leq&C'\(\|v\|^2+ (\lam-\lam_0)^2\)\hbox{ for some constant }C',\end{split}\end{equation} where $\p>0$ is the bound of  $u\in\cM_{loc}^\lam(\.)$, which is independent of $\lam$ by \eqref{eeq2}.
Taking the inner product of the equation \eqref{eq8} with $v$ and using \eqref{eq1.4} to obtain that there is a constant $C>0$ being independent of $\lam$ such that \begin{equation*}\begin{split}\frac{d}{dt}\|v\|^2
 &\leq C\(\|v\|^2+ (\lam-\lam_0)^2\).\end{split}\end{equation*}
Applying  the classical Gronwall lemma to get that
$$\|v(t)\|^2\leq \(e^{Ct}-1\)(\lam-\lam_0)^2,$$
which completes the proof. $\Box$


\bl\label{le1} Under the assumptions {\bf (F1),\,(F2)}, $y=0$ is locally asymptotically stable for $\phi_{\lam_0}$. Therefore $0$ is a pullback attractor of $\phi_{\lam_0}$. \el
{\bf Proof.} Since $X_c^\a$ is finite dimensional,  all the norms on $X_c^\a$ are equivalent. Hence for convenience, we equip $X^\a_c$ the norm $\|\.\|$ of $X$ in the following argument.

Note that $\phi_{\lam_0}$ is generated by the equation \be\label{e2}y_t=p(t,y+\xi^{\lam_0}_{\0_tp}(y)),\hs y\in \overline{\mb{B}}_{X^\a_c}(\varrho),\,p\in \cH.\ee Taking the inner product of the equation
\eqref{e2} with $y+\xi^{\lam_0}_{\0_tp}(y)$ in $X$, using the fact $(y,\,\xi^{\lam_0}_{\0_tp}(y))=0$ and the assumption {\bf(F1)}, it yields
\begin{equation}\label{eq1.3}\begin{split}\frac{1}{2}\frac{d}{dt}\|y\|^2&=\(p\(t,y+\xi^{\lam_0}_{\0_tp}(y)\),\,y+\xi^{\lam_0}_{\0_tp}(y)\)\\
&\leq-\b\.\kappa\(y+\xi^{\lam_0}_{\0_tp}(y)\).\end{split}\end{equation}
It is clear that $\kappa\(y+\xi^{\lam_0}_{\0_tp}(y)\)=0$ if and only if $y=0$. Therefore $\lim_{t\ra\8}\|y\|=0.$ The proof is complete. $\Box$

Henceforth we will suppose that \benu\item[{\bf (F3)}] The hull $\cH$ is a compact  metric space. \eenu

We then obtain a pullback attractor bifurcation theory for $\phi_\lam$ as $\lam$ crosses $\lam_0$.
\bt\label{le2} Under the assumptions {\bf (F1),\,(F2)} and {\bf(F3)}, the cocycle semiflow $\phi_\lam$ bifurcates from $(0,\lam_0)$ a pullback attractor $A_\lam(\.)$ for $\lam>\lam_0$, and
\be\label{eq12}\lim_{\lam\ra\lam_0^+}H_{X_c^\a}(A_\lam(p),\{0\})=0.\ee
\et
{\bf Proof.}  
Recall from Lemma \ref{le1} that $0$ is a pullback attractor for $\phi_{\lam_0}$ and it pullback attracts $\overline{\mb{B}}_{X^\a_c}(\varrho)$ for sufficiently small $\varrho>0$. The bounded set $\overline{\mb{B}}_{X^\a_c}(\varrho)\<X^\a_c$ is compact due to $X^\a_c$ being  finite dimensional.  Moreover, $\overline{\mb{B}}_{X^\a_c}(\varrho)$  is forward invariant under $\phi_{\lam_0}$. Then by Corollary \ref{c1}, there is a $\eta'\in(0,\eta/8)$ such that for each $\lam\in I_{\lam_0}(\eta')$, the cocycle semiflow $\phi_\lam$ has a pullback attractor $A_\lam(\.)$ and \eqref{eq12} holds.

In the following, we prove that $0\notin A_\lam(\.)$ for  $\lam\in I^+_{\lam_0}(\eta'):=(\lam_0,\lam_0+\eta')$, which completes the proof.

Let $\lam\in I^+_{\lam_0}(\eta')$ be fixed, and let $w(t)=y(-t)$. Then $w(t)$ satisfies \be\label{e5}w_t-(\lam_0-\lam)w=-p(-t,w+\xi^\lam_{\0_{-t}p}(w)).\ee Taking the inner product of the equation
\eqref{e5} with $w$ in $X^\a$, we have \be\label{e6}\frac{1}{2}\frac{d}{dt}\|w\|^2-(\lam_0-\lam)\|w\|^2=-(p(t,w+\xi^\lam_{\0_{-t}p}(w)),w).\ee
Since $$\|p(t,u)\|\leq k(\|u\|_\a)\|u\|_\a\hbox{ and }\|\xi^\lam_{\0_{-t}p}(w)\|_\a\leq L_1\|w\|_\a,$$
we have\begin{equation}\label{e7}\begin{split}\|p(-t,w+\xi^\lam_{\0_{-t}p}(w))\|\leq& k(\|w+\xi^\lam_{\0_{-t}p}(w)\|_\a)\|w+\xi^\lam_{\0_{-t}p}(w)\|_\a\\
\leq& k(\|w+\xi^\lam_{\0_{-t}p}(w)\|_\a)\(\|w\|_\a+L_1\|w\|_\a\)\\
\leq&k(\|w+\xi^\lam_{\0_{-t}p}(w)\|_\a)\.(1+L_1)\|w\|_\a\\
\leq&[(1+L_1)ck(\|w+\xi^\lam_{\0_{-t}p}(w)\|_\a)]\.\|w\|\\
\leq&\frac{1}{2}(\lam-\lam_0)\|w\|,\hbox{ for sufficiently small }\|w\|_\a,\end{split}
\end{equation}where $c>0$ only depends on $\a$. We get from \eqref{e6} and \eqref{e7} that $$\frac{d}{dt}\|w\|^2\leq-(\lam-\lam_0)\|w\|^2$$ for sufficiently small $\|w\|_\a$, which shows for fixed $\lam\in I_{\lam_0}^+(\eta')$, $0$ locally asymptotically stable for the cocycle flow generated by the equation \eqref{e5}. In other words, $0$ is a repeller of $\phi_\lam$ when $\lam\in I^+_{\lam_0}(\eta')$ and repels a neighborhood of $0$ in $X^\a_c$. This implies that $0\notin A_\lam(\.)$, $\lam\in I^+_{\lam_0}(\eta')$. The proof is complete. $\Box$

We are now in position to give and prove the main result of this paper concerning the invariant set bifurcation of $\va_\lam$.
\bt\label{th2}  Under the assumptions {\bf (F1),\,(F2)} and {\bf(F3)}, the cocycle semiflow $\va_\lam$ bifurcates from $(0,\lam_0)$ an invariant compact set $B_\lam(\.)$ for $\lam>\lam_0$, and for each $p\in P$,\be\label{e11}\lim_{\lam\ra\lam_0^+}H_X(B_\lam(p),\{0\})=0.\ee\et
{\bf Proof.} Let $A_\lam(\.)$ be the bifurcated attractor obtained in Theorem \ref{le2}. Define $B_\lam(\.)$ by \be\label{eq1.8}B_\lam(p)=\{y+\xi^\lam_p(y):y\in A_\lam(p)\},\hs p\in\cH.\ee
%

We know from Theorem \ref{le2} that $0\notin B_\lam(\.)$ and $B_\lam(\.)\<\cM_{loc}^\lam(\.)$.
Based on the compactness of  $A_\lam(p)$ and the continuity of $\xi^\lam_p(y)$ in $y$,  we can directly derive the compactness of $B_\lam(p)$. So $B_\lam(\.)$ is compact.

We claim that $B_\lam(\.)$ is invariant under $\va_\lam$. Indeed,  let $p\in P$ and $y+\xi^\lam_p(y)\in B_\lam(p)$.
Since $\phi_\lam(t,p)y\in A_\lam(\0_tp)$, $t\geq0$, by the invariance of $\cM_{loc}^\lam(\.)$, we have
$$\va_\lam(t,p)(y+\xi^\lam_p(y))=\phi_\lam(t,p)y+\xi^\lam_{\0_tp}(\phi_\lam(t,p)y)\in B_\lam(\0_tp),$$
which shows $$\va_\lam(t,p)B_\lam(p)\<B_\lam(\0_tp),\hs t\geq0.$$
On the other hand, for any $y+\xi^\lam_{\0_tp}(y)\in  B_\lam(\0_tp)$, $t\geq0$. Using the invariance of $A_\lam(\.)$ and  $\cM_{loc}^\lam(\.)$, there is a $y'\in A_\lam(p)$ such that $y=\phi_\lam(t,p)y'$. Then \begin{equation*}\begin{split}y+\xi^\lam_{\0_tp}(y)&=\phi_\lam(t,p)y'+\xi^\lam_{\0_tp}(\phi_\lam(t,p)y')\\
&=\va_\lam(t,p)(y'+\xi^\lam_{\0_tp}(y'))\in \va(t,p)B_\lam(p),\end{split}\end{equation*}
which shows $$B_\lam(\0_tp)\<\va_\lam(t,p)B_\lam(p),\hs t\geq0.$$Therefore $B_\lam(\.)$ is invariant under $\va_\lam$.

Finally, \eqref{e11} is an immediately consequence of \eqref{eq12} and \eqref{eq1.5}. $\Box$

We now give a result which  parallels Theorem \ref{th2}.
\bc  Let the assumptions {\bf (F1),(F2),(F3)} hold, but replace \eqref{eq1.7} by the
assumption that $$\(f(t,u),u\)\geq\b\.\kappa(u).$$ Then the cocycle semiflow $\va_\lam$ bifurcates from $(0,\lam_0)$ an invariant compact set $B_\lam(\.)$ for $\lam<\lam_0$, and for each $p\in P$,\be\label{e11'}\lim_{\lam\ra\lam_0^-}H_X(B_\lam(p),\{0\})=0.\ee\ec
{\bf Proof.} Let $\lam\in I_{\lam_0}(\eta/8)$. Consider the following equation \be\label{eq3'}z_t-(\lam_0-\lam)z=-p(-t,z+\xi^\lam_{\0_{-t}p}(z)),\hs z\in \overline{\mb{B}}_{X^\a_c}(\varrho),\,p\in\cH.\ee Denote by $\phi^-_\lam$ be the cocycle flow    generated by \eqref{eq3'}. Then $\phi^-_\lam$  be the inverse flow of $\phi_\lam$.

Repeating the argument of Lemma \ref{le3}, Lemma \ref{le1} and Theorem \ref{le2} (replacing $\phi_\lam$  by $\phi_\lam^-$) to show $\phi_\lam^-$ bifurcates from $(0,\lam_0)$ a pullback attractor $R_\lam(\.)$ for $\lam<\lam_0$, and
\be\label{eq12}\lim_{\lam\ra\lam_0^-}H_{X_c^\a}(R_\lam(p),\{0\})=0.\ee
It is clear that $R_\lam(\.)$ is also an invariant set of $\phi_\lam$. Define a set $B_\lam(\.)$ by $$B_\lam(p)=\{y+\xi^\lam_p(y):y\in R_\lam(p)\},\hs p\in\cH.$$
Similar to Theorem \ref{th2}, we can show $B_\lam(\.)$ is an invariant set of $\va_\lam$ and \eqref{e11'} holds. $\Box$
\section{An example}
Consider the nonautonomous system:
\be\label{eq7.1}\left\{\ba{lll}u_t-\Delta u=\lam u\pm h(t)u^3,\hs \,t>0,x\in\W;\\[1ex]\Hs\hs u=0,\Hs\Hs\Hs t>0, x\in\partial\W,\ea\right.\ee where $\W$ is a bounded domain in $\R^3$ with smooth boundary, $h$ is a function such that $h(t)\geq\de>0$ for some $\de>0$.

Denote by $A$ the operator $-\Delta$ associated with the homogeneous  Dirichlet boundary condition.  Then $A$ is a sectorial operator on $X=L^2(\Omega)$ with compact resolvent, and $\mathcal {D}(A)=H^2(\W)\bigcap H_0^1(\W)$.
Note that $A$  has eigenvalues
$0<\mu_1<\mu_2<\cdots<\mu_k<\cdots$. Denote $V=H_0^1(\Omega)$. By $(\cdot,\cdot)$ and
$|\cdot|$ we denote the usual inner product and norm on $H$,
respectively. The inner product and norm on $V$, denoted by $((\.,\.))$ and $\|\cdot\|$, respectively, are defined as
$$
((u,v))= \int_{\Omega}\nabla u\.\nab v\mathrm{d}x,\hs  \|u\|=\(\int_{\Omega}|\nabla u|^2\mathrm{d}x\)^{1/2}
$$for $u,v\in V$.

The system $(\ref{eq7.1})$ can be written into  an abstract  equation on $X$:  $$u_t+Au=\lam u\pm h(t)u^3.$$
Define  the hull $\cH:=\cH[h(\.)u^3]$.
By the assumption on $h$, it is clear that $$(p(t,u),u)\geq\de\int_\W u^4dx,\hs p\in \cH.$$
Consider the cocycle system:
\be\label{eq7.2}
u_t+Au=\lam u \pm p(t,u), \Hs p\in\cH.\ee Denote $\va^\pm_\lam:=\va^\pm_\lam(t,p)u$
 the cocycle semiflow on $H^1_0(\Omega)$  driven by the base flow (translation group) $\0$ on  $\cH$.

Since all the hypotheses in the main theorem above  are fulfilled,  we     obtain some interesting results concerning the dynamics of the perturbed system. In particular,
\bt Suppose $\cH$ is compact. Then the cocycle semiflow $\va^-_\lam$ (resp. $\va^+_\lam$) bifurcates from $(0,\mu_k)$, $k=1,2,\cdots$ an invariant compact set $B^-_\lam(\.)$ for $\lam>\lam_0$ (resp. $B^+_\lam(\.)$ for $\lam<\lam_0$) and for each $p\in P$,$$\lim_{\lam\ra\lam^+_0}H_X(B^-_\lam(p),\{0\})=0.$$
$$\big(\hbox{resp. }\lim_{\lam\ra\lam^-_0}H_X(B^+_\lam(p),\{0\})=0.\big)$$\et

\section{Appendixes}

\subsection{Relationship between the pullback attraction of $\va$ and the attraction of $\Phi$}
\noindent{\bf Proof of Proposition \ref{P1}.}
{\bf Necessity:} By the compactness of $P$, one finds that \begin{equation*}\begin{split}\lim_{t\ra\8}H_Y\(\Phi(t)\mathbb{B},P\X K_P\)=&\lim_{t\ra\8}H_X\(\va(t,p)B(p),K_P\)\\
\leq&\lim_{t\ra\8}\sup_{p\in P}H_X\(\va(t,p)B(p),K_P\)\\
=&\lim_{t\ra\8}\sup_{p\in P}H_X\(\va(t,\0_{-t}p)B(\0_{-t}p),K_P\)\\
=&0.\\
\end{split}
\end{equation*} This means the compact set $P\X K_P$ attracts $\mathbb{B}$ through $\Phi$. Therefore the omega-limit
set $\w(\mathbb{B})$ of $\mathbb{B}$ exists and attracts $\mathbb{B}$.

In the following, we prove $\w(\mathbb{B})\<\mathbb{K}$, which completes the necessity.
For this purpose, define a nonautonomous set $\tilde{B}(\.)$ as follows $$\tilde{B}(p):=\overline{\bigcup_{s\geq0}\va(s,\0_{-s}p)B(\0_{-s}p)},\hs p\in P.$$
It is clear that $B(\.)\<\tilde{B}(\.)$. We first say $\tilde{B}(\.)$ is forward invariant. Indeed, for any $t\geq0$ and $p\in P$,
\begin{equation}\label{eq2.5}
\begin{split}
\va(t,p)\tilde{B}(p)&=\va(t,p)\overline{\bigcup_{s\geq0}\va(s,\0_{-s}p)B(\0_{-s}p)}\\
&\<\overline{\bigcup_{s\geq0}\va(t,p)\circ\va(s,\0_{-s}p)B(\0_{-s}p)}\\
&=\overline{\bigcup_{s\geq0}\va(t+s,\0_{-(t+s)}\circ\0_tp)B(\0_{-(t+s)}\circ\0_tp)}\\
&\<\overline{\bigcup_{s\geq0}\va(s,\0_{-s}\circ\0_tp)B(\0_{-s}\circ\0_tp)}=\tilde{B}(\0_tp).\end{split}
\end{equation} So $\tilde{B}(\.)$ is forward invariant, which implies the omega-limit
set $\w({\tilde{B}})(\.)$ of $\tilde{B}(\.)$ is the maximal invariant set in $\tilde{B}(\.)$.
Furthermore, for each $p\in P$, \begin{equation*}
\begin{split}\w(\tilde{B})(p)&=\bigcap_{\tau\geq0}\overline{\bigcup_{t\geq\tau}\va(t,\0_{-t}p)\tilde{B}(\0_{-t}p)}\\
&=\bigcap_{\tau\geq0}\overline{\bigcup_{t\geq\tau}\va(t,\0_{-t}p)\overline{\bigcup_{s\geq0}\va(s,\0_{-(s+t)}p)B(\0_{-(s+t)}p)}}\\
&=\bigcap_{\tau\geq0}\overline{\bigcup_{t\geq\tau}\va(t,\0_{-t}p)\bigcup_{s\geq0}\va(s,\0_{-(s+t)}p)B(\0_{-(s+t)}p)}\\
&=\bigcap_{\tau\geq0}\overline{\bigcup_{t\geq\tau}\bigcup_{s\geq0}\va(t,\0_{-t}p)\circ\va(s,\0_{-(s+t)}p)B(\0_{-(s+t)}p)}\\
&=\bigcap_{\tau\geq0}\overline{\bigcup_{t\geq\tau}\bigcup_{s\geq0}\va(t+s,\0_{-(s+t)}p)B(\0_{-(s+t)}p)}\\
&=\bigcap_{\tau\geq0}\overline{\bigcup_{t\geq\tau}\va(t,\0_{-t}p)B(\0_{-t}p)}\\
&=\w(B)(p),\end{split}
\end{equation*}
where the third ``$=$'' holds since for each fixed $t\geq0$ and $p\in P$, $\va(t,\0_{-t}p)$ is a continuous map on $X$.
It follows that $\w(B)(\.)$ is the maximal forward invariant set in $\tilde{B}(\.)$. Therefore $\mathbb{C}:=\bigcup_{p\in P}\big(\{p\}\X \w(B)(p)\big)$  is the maximal invariant set in $\mathbb{\tilde{B}}:=\bigcup_{p\in P}\big(\{p\}\X \tilde{B}(p)\big)$.
By the forward invariance of $\tilde{B}(\.)$, \begin{equation*}
\begin{split}\va(t)\mathbb{\tilde{B}}&=\va(t)\bigcup_{p\in P}\big(\{p\}\X \tilde{B}(p)\big)\\
&\<\bigcup_{p\in P}\va(t)\big(\{p\}\X \tilde{B}(p)\big)\\
&=\bigcup_{p\in P}\big(\{\0_tp\}\X \va(t,p)\tilde{B}(p)\big)\\
&\<\hbox{(by \eqref{eq2.5})}\<\bigcup_{p\in P}\big(\{\0_tp\}\X \tilde{B}(\0_tp)\big)\\
&=\mathbb{\tilde{B}},\hs t\geq0,\end{split}
\end{equation*}
i.e. $\mathbb{\tilde{B}}$ is positively invariant under $\va$. Then $\w(\mathbb{\tilde{B}})$ is also the  maximal invariant set in $\mathbb{\tilde{B}}$.
Therefore  \be\label{eq2.2}\w(\mathbb{B})\<\w(\mathbb{\tilde{B}})=\mathbb{C}.\ee
Finally, by the assumption that $K(\.)$ attracts $B(\.)$, one knows that $\w(B)(\.)\<K(\.)$, and thus $\mathbb{C}\<\mathbb{K}$, which shows $$\w(\mathbb{B})\<\mathbb{K}.$$

\Vs
\noindent {\bf Sufficiency:} In a very similar way as above, we can prove the sufficiency.

By the compactness of $P$, \begin{equation*}\begin{split}
\lim_{t\ra\8}H_X\(\va(t,\0_{-t}p)B(\0_{-t}p),K_P\)]\leq&\lim_{t\ra\8}\sup_{p\in P}H_X\(\va(t,p)B(p),K_P\)\\
=&\lim_{t\ra\8}\sup_{p\in P}H_Y\(\Phi(t)\mathbb{B},P\X K_P\)\\
=&\lim_{t\ra\8}H_Y\(\Phi(t)\mathbb{B},P\X K_P\)\\
=&0,\\
\end{split}
\end{equation*} which implies $\w(B)(\.)$ exists and pullback attracts $B(\.)$.

To complete the proof, it suffices to show $\w(B)(\.)\<K(\.)$.
We first define a set $$\hat{\mathbb{B}}=\overline{\bigcup_{s\geq0}\Phi(s)\mathbb{B}}.$$ Then $\hat{\mathbb{B}}$ is positively invariant and $$\w(\hat{\mathbb{B}})=\w(\mathbb{B}).$$ This implies that $\w(\mathbb{B})$ is the  maximal invariant set in $\hat{\mathbb{B}}$.
Write $\w(\mathbb{B}):=\bigcup_{p\in P}\{p\}\times C(p)$, then $C(\.)$ is the  maximal invariant set in $\hat{B}(\.)$, where $\hat{B}(\.)$ is the set defined by $\hat{\mathbb{B}}:=\bigcup_{p\in P}\{p\}\times \hat{B}(p)$. By the positive invariance of $\hat{\mathbb{B}}$, one also knows that $\hat{B}(\.)$ is forward invariant. This implies $\w({\hat{B}})(\.)$ is the  maximal invariant set in $\hat{B}(\.)$. We then have that $\w(B)(\.)\<\w(\hat{B})(\.)=C(\.)$. We learn from the condition $\w(\mathbb{B})\<\mathbb{K}$ that $C(\.)\<K(\.)$. In summary, $\w(B)(\.)\<K(\.)$, which completes the sufficiency.
$\Box$
\subsection{Construction of local invariant manifold}
Let $M>0$. For $\mu\geq 0$, define a Banach space as
$$\mathscr{X}_\mu=\left\{u\in C(\R;X^\a): \,\,\sup_{t\in\R}e^{-\mu |t|}\|x(t)\|_\a\leq M\right\},$$ which is equipped with the norm $\|\.\|_{\mathscr{X}_\mu}$, $$\|x\|_{\mathscr{X}_\mu}=\sup_{t\in\R}e^{-\mu |t|}\|x(t)\|_\a,\Hs \A\,x\in \sX_\mu^\a.$$

Let $A^\lam=A-\lam$. Write $\sig(A_\lam)=\sig_-\cup \sig_c\cup\sig_+$, where $$\sig_c=\{\lam_0-\lam\},$$
$$\sig_-=\sig(A_\lam)\cap\{\hbox{Re}\, \lambda<0\},\hs \sig_+=\sig(A_\lam)\cap\{\hbox{Re}\, \lambda>0\}.$$
According to the  spectral decomposition, the space $X$ has a direct sum
decomposition: $X=X_-\oplus X_c\oplus X_+$. Denote $X_\pm:=X_-\bigcup X_+$. Note that each $X_i$, $i=-,+,\pm,c$ is independent of $\lam$.
Let $$\Pi_i:X\ra X_i,\Hs i=-,+,\pm,c$$ be the projection from $X$ to $X_i$. Denote $A^\lam_i=A^\lam|_{X_i}$. By the assumption {\bf(F2)}, we deduce that if $\lam\in(\lam_0-\eta/4,\lam_0+\eta/4)$ then for $\a\in[0,1)$,
\be\label{b1*}\|A^\a e^{-A^\lam_-t}\|\leq e^{\frac{3\eta}{4}t},\hs \|e^{-A^\lam_-t}\|\leq e^{-\frac{3\eta}{4}t},\Hs t\leq0,\ee
\be\label{b2*}\|A^\a e^{-A^\lam_+t}\Pi_+A^{-\a}\|\leq e^{-\frac{3\eta}{4}t},\hs\|A^\a e^{-A^\lam_+t}\|\leq t^{-\a}e^{-\frac{3\eta}{4}t},\Hs t>0,\ee
\be\label{b3*}\|A^\a e^{-A^\lam_ct}\|\leq e^{\frac{\eta}{4}|t|},\hs\| e^{-A^\lam_ct}\|\leq e^{\frac{\eta}{4}|t|},\Hs t\in\R.\ee

\Vs{\bf Proof of Proposition \ref{p2}.} Let $\chi:\R\ra\R$ be a smooth function such that $$\chi(z)=\left\{\ba{ll}
  1, \Hs\Hs |z|\leq1/2;\\[1ex]
  0, \Hs\Hs |z|\geq1.\ea\right. $$
For $\p>0$, one can then define a smooth function
such that $$p_\p(t,u)=\chi\(\frac{\|u\|_\a}{\p}\)p(t,u).$$Select suitable $\chi$ such that  \be\label{eq2.4}\|p_\p(t,u)-p_\p(t,v)\|\leq k(\p)\|u-v\|,\ee where $k(\p)$ is the local Lipschitz constant of $f$ given in \eqref{eqq6}.
Instead of \eqref{e1}, we consider the truncated system
\be\label{eq1.1}u_t+Au=\lam u+p_\p(t,u),\Hs p\in\cH.\ee
Suppose that $\p$ is so small that \be\label{eq1.6}M_\p:=k(\p)\int_0^\8\(2+\tau^{-\a}\) e^{-\frac{\eta}{4} \tau}d\tau<1.\ee

Let  $u\in\mathscr{X}_{\eta/2}$. By simple computations, we know  that $u$ is the solution of \eqref{eq1.1} if and only if it solves the integral equation
\begin{equation}\label{e4}
\begin{split}
u(t)&=e^{-A^\lam_ct}\Pi_cu(0)+\int_{0}^te^{-A^\lam_c(t-\tau)}\Pi_cp_\p(\tau,u(\tau))d\tau\\
&\quad+\int_{-\8}^{t}e^{-A^\lam_+(t-\tau)}\Pi_+p_\p(\tau,u(\tau))d\tau\\
&\quad-\int_{t}^\8e^{-A^\lam_-(t-\tau)}\Pi_-p_\p(\tau,u(\tau))d\tau.
\end{split}
\end{equation}
Take a $\tilde{\varrho}>0$ small enough so that
\be\label{equ2}\tilde{\varrho}\leq\(1-M_\p\)M.\ee  Let $p\in\cH$ and $\lam\in I_{\lam_0}(\eta/8)$ be fixed. For each $y\in \overline{\mb{B}}_{X^\a_c}(\tilde{\varrho})$, one can use the righthand side of equation $(\ref{e4})$ to define a contraction mapping $\cT:=\cT_y$ on
 $\mathscr{X}_{\eta/2}$ as follows:
\begin{equation*}
\begin{split}
\cT u(t)&=e^{-A^\lam_ct}y+\int_{0}^te^{-A^\lam_c(t-\tau)}\Pi_cp_\p(\tau,u(\tau))d\tau\\
&\quad+\int_{-\8}^{t}e^{-A^\lam_+(t-\tau)}\Pi_+p_\p(\tau,u(\tau))d\tau\\
&\quad-\int_{t}^\8e^{-A^\lam_-(t-\tau)}\Pi_-p_\p(\tau,u(\tau))d\tau.
\end{split}
\end{equation*}
We first verify that $\mathcal{T}$ maps $\mathscr{X}_{\eta/2}$ into itself.

For notational convenience, we write $$0\wedge t=\min\{0,t\},\hs 0\vee t=\max\{0,t\},\hs \hbox{for }t\in\R.$$
Let $u\in \mathscr{X}_{\eta/2}$. By \eqref{b1*}-\eqref{b3*} and \eqref{eq2.4} we have
 \begin{equation}\label{eq1.9}
\begin{split}
\|\mathcal {T}u(t)\|_\a&\leq e^{\frac{\eta}{4}|t|}\|y\|_\a+\int_{0\wedge t}^{0\vee t}e^{\frac{\eta}{4}|t-\tau|}k(\p)\|u(\tau)\|_\a d\tau\\
&\quad+ \int_{-\8}^t (t-\tau)^{-\a} e^{-\frac{3\eta}{4}(t-\tau)}k(\p)\|u(\tau)\|_\a d\tau\\
&\quad+\int_t^\8  e^{\frac{3\eta}{4}(t-\tau)}k(\p)\|u(\tau)\|_\a d\tau.
\end{split}\end{equation}
It is trivial to verify that \begin{equation}\label{eq2.0}\begin{split}
&\hs e^{-\frac{\eta}{2}|t|}\int_{0\wedge t}^{0\vee t}e^{\frac{\eta}{4}|t-\tau|}k(\p)\|u(\tau)\|_\a d\tau\\
&=\int_{0\wedge t}^{0\vee t}e^{-\frac{\eta}{4}|t-\tau|}\big[e^{-\frac{\eta}{2}|\tau|}k(\p)\|u(\tau)\|_\a\big] d\tau.\end{split}
\end{equation}
Observing that
\begin{equation*}
\begin{split}
e^{-\frac{\eta}{2}|t|}=e^{-\frac{\eta}{2}|(t-\tau)+\tau|}\leq e^{\frac{\eta}{2}|t-\tau|}e^{-\frac{\eta}{2}|\tau|},
\end{split}
\end{equation*}
by \eqref{equ2}, \eqref{eq1.9} and \eqref{eq2.0} we find that
\begin{equation}\label{eq3.1}
\begin{split}
&\hs e^{-\frac{\eta}{2}|t|}\|\mathcal {T}x(t)\|_{\a}\leq e^{-\frac{\eta}{4}|t|}\|y\|_{\a}\\[1ex]
&\quad+\int_{0\wedge t}^{0\vee t}e^{-\frac{\eta}{4}|t-\tau|}\big[e^{-\frac{\eta}{2}|\tau|}k(\p)\|u(\tau)\|_\a\big] d\tau\\[1ex]
&\quad+ \int_{-\8}^t (t-\tau)^{-\a}e^{\frac{\eta}{2}|t-\tau|} e^{-\frac{3\eta}{4}(t-\tau)}\big[e^{-\frac{\eta}{2}|\tau|}k(\p)\|u(\tau)\|_\a\big]d\tau\\[1ex]
&\quad+\int_t^\8 e^{\frac{\eta}{2}|t-\tau|}e^{\frac{3\eta}{4}(t-\tau)}\big[e^{-\frac{\eta}{2}|\tau|}k(\p)\|u(\tau)\|_\a\big] d\tau\\[1ex]
&=e^{-\frac{\eta}{4}|t|}\|y\|_{\a}+\int_{0\wedge t}^{0\vee t}e^{-\frac{\eta}{4}|t-\tau|}\big[e^{-\frac{\eta}{2}|\tau|}k(\p)\|u(\tau)\|_\a\big]d\tau\\[1ex]
&\quad+\int_{-\8}^t (t-\tau)^{-\a} e^{-\frac{\eta}{4}(t-\tau)}\big[e^{-\frac{\eta}{2}|\tau|}k(\p)\|u(\tau)\|_\a\big]d\tau\\[1ex]
&\quad+\int_t^\8  e^{\frac{\eta}{4}(t-\tau)}\big[e^{-\frac{\eta}{2}|\tau|}k(\p)\|u(\tau)\|_\a\big] d\tau\\[1ex]
&\leq \|y\|_\a+M_\p\|u\|_{\mathscr{X}_{\eta/2}}\leq M,\hs \A\,t\in\R.
\end{split}
\end{equation}
 Hence   $\mathcal{T}u\in\mathscr{X}_{\eta/2}$.

Next, we check that $\mathcal{T}$ is contractive. Indeed,  in a quite similar fashion as above,  it can be shown that for any $u,u'\in\mathscr{X}_{\eta/2}$,
\begin{equation}\label{eq2.2}
\begin{split}
&\hs\,\, e^{-\frac{\eta}{2}|t|}\|\mathcal {T}u(t)-\mathcal {T}u'(t)\|_\a\\
&\leq k(\p)\int_{0\wedge t}^{0\vee t}e^{-\frac{\eta}{4}|t-\tau|}\(e^{-\frac{\eta}{2}|\tau|}\|u(\tau)-u'(\tau)\|_\a \)d\tau\\
&\quad+ k(\p)\int_{-\8}^t(t-\tau)^{-\a} e^{-\frac{\eta}{4}(t-\tau)}\(e^{-\frac{\eta}{2}|\tau|}\|u(\tau)-u'(\tau)\|_\a \)d\tau\\
&\quad+k(\p)\int^\8_t e^{\frac{\eta}{4}(t-\tau)}\(e^{-\frac{\eta}{2}|\tau|}\|u(\tau)-u'(\tau)\|_\a\)d\tau\\
&\leq \(k(\p)\int_0^\8\(2+\tau^{-\a}\) e^{-\frac{\eta}{4} \tau}d\tau\)\|u-u'\|_{\mathscr{U}_{\eta/2}}\\
&:=M_\p\|u-u'\|_{\mathscr{X}_{\eta/2}},\Hs\A\, t\in\R.
\end{split}
\end{equation}
Thus
$$
\|\mathcal {T}u-\mathcal {T}u'\|_{\mathscr{X}_{\eta/2}}\leq M_\p\|u-u'\|_{\mathscr{X}_{\eta/2}}.
$$
The conditon \eqref{eq1.6} then asserts that $\cT$ is contractive.

Thanks to the Banach fixed-point theorem, $\cT$ has a unique fixed point $\gamma^y_{p,\lam}\in \mathscr{X}_{\eta/2}$ which   is precisely a full solution of $(\ref{e1})$ with $\Pi_c\gamma^y_{p,\lam}(0)=y$ and solves the integral equation
\begin{equation}\label{equ5}
\begin{split}
\gamma^y_{p,\lam}(t)&=e^{-A^\lam_ct}y+\int_{0}^te^{-A^\lam_c(t-\tau)}\Pi_cp_\p(\tau,\gamma^y_{p,\lam}(\tau))d\tau\\
&\quad+\int_{-\8}^{t}e^{-A^\lam_+(t-\tau)}\Pi_+p_\p(\tau,\gamma^y_{p,\lam}(\tau))d\tau\\
&\quad-\int_{t}^\8e^{-A^\lam_-(t-\tau)}\Pi_-p_\p(\tau,\gamma^y_{p,\lam}(\tau))d\tau.
\end{split}
\end{equation}

For    $y,z\in \overline{\mb{B}}_{X^\a_c}(\tilde{\varrho})$ and $t\in\R$, similar to \eqref{eq2.2}, by  $(\ref{equ5})$ we have
\begin{equation*}
\begin{split}
e^{-\frac{\eta}{2} |t|}&\|\gamma^y_{p,\lam}(t)-\gamma^z_{p,\lam}(t)\|_\a\leq\,e^{-\frac{\eta}{4} |t|}\|y-z\|_\a\\
&+ k(\p)\int_{0\wedge t}^{0\vee t}e^{-\frac{\eta}{4}|t-\tau|}\big(e^{-\frac{\eta}{2}|\tau|}\|\gamma^y_{p,\lam}(\tau)-\gamma^z_{p,\lam}(\tau)\|_\a \big)d\tau\\
&+ k(\p)\int_{-\8}^t(t-\tau)^{-\a} e^{-\frac{\eta}{4}(t-\tau)}\big(e^{-\frac{\eta}{2}|\tau|}\|\gamma^y_{p,\lam}(\tau)-\gamma^z_{p,\lam}(\tau)\|_\a \big)d\tau\\
&+k(\p)\int^\8_t e^{\frac{\eta}{4}(t-\tau)}\big(e^{-\frac{\eta}{2}|\tau|}\|\gamma^y_{p,\lam}(\tau)-\gamma^z_{p,\lam}(\tau)\|_\a \big)d\tau\\
\leq\,& \|y-z\|_\a+ M_\p\|\gamma^y_{p,\lam}-\gamma^z_{p,\lam}\|_{\mathscr{X}_{\eta/2}}.
\end{split}
\end{equation*}
Hence $$
\|\gamma^y_{p,\lam}-\gamma^z_{p,\lam}\|_{\mathscr{X}_{\eta/2}}\leq \frac{1}{1-M_\p}\|y-z\|_\a,
$$
which implies that \be\label{eq2}\|\gamma^y_{p,\lam}(0)-\gamma^z_{p,\lam}(0)\|_\a\leq  \frac{1}{1-M_\p}\|y-z\|_\a.\ee

For each $p\in\cH$ and $\lam\in I_{\lam_0}(\eta/8)$, define a mapping from $X_c^\a$ to $X_{us}^\a$ as \begin{equation}
\begin{split}\label{eq2.5}
\xi^\lam_p(y):=&\int_{-\8}^0e^{A^\lam_+\tau}\Pi_+p_\p(\tau,\gamma^y_{p,\lam}(\tau))d\tau\\
&-\int_0^\8e^{A^\lam_-\tau}\Pi_-p_\p(\tau,\gamma^y_{p,\lam}(\tau))d\tau,\hs y\in\overline{\mb{B}}_{X^\a_c}(\tilde{\varrho}).
\end{split}
\end{equation}
Setting $t=0$ in \eqref{equ5} leads to
\be\label{eq1.2}\gamma^y_{p,\lam}(0)=y+\xi^\lam_p(y),\Hs y\in \overline{\mb{B}}_{X^\a_c}(\tilde{\varrho}).\ee
We conclude from \eqref{eq2}, \eqref{eq2.5} and \eqref{eq1.2} that $\xi^\lam_p(\.):\overline{\mb{B}}_{X^\a_c}(\tilde{\varrho})\ra X^\a_{us}$ is a Lipschitz continuous mapping uniformly on $p$ and $\lam$. More specifically, let $$L_1:=\frac{1}{1-M_\p}+1.$$Then for each $y,z\in \overline{\mb{B}}_{X^\a_c}(\tilde{\varrho})$, \begin{equation*}
\begin{split}\|\xi^\lam_p(y)-\xi^\lam_p(z)\|_\a&\leq\|\gamma^y_{p,\lam}(0)-\gamma^z_{p,\lam}(0)\|_\a+\|y-z\|_\a\\
&\leq L_1\|y-z\|_\a.\end{split}
\end{equation*}
Since $\gamma^y_{p,\lam}\equiv0$ is a full solution of \eqref{eq1.1}, we have $\xi_p^\lam(0)\equiv0$, and thus $$\lim_{\|y\|_\a\ra0}\|\xi^\lam_p(y)\|_\a=0$$ uniformly on   $p\in\cH$ and $\lam\in I_{\lam_0}(\eta/8)$.

Take a sufficiently small $\varrho>0$ such that $$\|y+\xi_p^\lam(y)\|\leq \frac{\p}{2},\hs y\in \overline{\mb{B}}_{X^\a_c}(\varrho).$$ Define for each $p\in\cH$ the $p$-section as $$\cM_{loc}^\lam(p)=\{y+\xi^\lam_p(y):y\in \overline{\mb{B}}_{X^\a_c}(\varrho)\}.$$  By the definition of $p_\p$, $\cM_{loc}^\lam(\.):=\{\cM_{loc}^\lam(p)\}_{p\in\cH}$ is a local invariant manifold of the cocycle semiflow $\va_\lam$, $\lam\in I_{\lam_0}(\eta/8)$ generated by \eqref{e1}. And for each $p\in\cH$, the section $\cM^\lam_{loc}(p)$ is homeomorphic to $\overline{\mb{B}}_{X^\a_c}(\varrho)$.

In the last part, we show $\xi_p^\.(y):I_{\lam_0}(\eta/8)\ra X^\a_{us}$ is Lipschitz uniformly on $y\in \overline{\mb{B}}_{X^\a_c}(\varrho)$ and $p\in P$. Indeed, for $\lam_1,\lam_2\in I_{\lam_0}(\eta/8)$ with $\lam_1\leq\lam_2$, we have for $t\in \R$ that \begin{equation*}
\begin{split}\|e^{-A^{\lam_1}_ct}-e^{-A^{\lam_2}_ct}\|&\leq \|e^{-A^{\lam_1}_ct}\|\.\big|1-e^{-(\lam_2-\lam_1)t}\big|\\
&\leq e^{\frac{\eta}{4}|t|}\.\big|1-e^{-(\lam_2-\lam_1)t}\big|.
\end{split}
\end{equation*}
Then for $t\in\R$,
\begin{equation}\label{eq5}
\begin{split}&e^{-\frac{\eta}{2} |t|}\int_{0}^t\|e^{-A^{\lam_1}_c(t-\tau)}p(\tau,\gamma^y_{p,\lam_1}(\tau))-e^{-A^{\lam_2}_c(t-\tau)}p_\p(\tau,\gamma^y_{p,\lam_2}(\tau))\|d\tau\\
\leq&\int_{0}^te^{-\frac{\eta}{4}|t-\tau|}k_1(\p)\(e^{-\frac{\eta}{2} |\tau|}\|\gamma^y_{p,\lam_1}(\tau)-\gamma^y_{p,\lam_2}(\tau)\|_\a\) d\tau\\
&+\int_{0}^te^{-\frac{\eta}{4}|t-\tau|}\.k_1(\p)\big|1 -e^{-(\lam_2-\lam_1)(t-\tau)}\big|\.\(e^{-\frac{\eta}{2}|\tau|}\|\gamma^y_{p,\lam_2}(\tau)\|_\a\) d\tau\\
\leq& k(\p)\int_{0}^te^{-\frac{\eta}{4}|t-\tau|}\(e^{-\frac{\eta}{2} |\tau|}\|\gamma^y_{p,\lam_1}(\tau)-\gamma^y_{p,\lam_2}(\tau)\|_\a\) d\tau\\
&+ k(\p) M\int_{0}^te^{-\frac{\eta}{4}|t-\tau|}\big|1 -e^{-(\lam_2-\lam_1)(t-\tau)}\big| d\tau.\end{split}
\end{equation}
We can apply very similar arguments as above to get that
\begin{equation}\label{eq6}
\begin{split}&e^{-\frac{\eta}{2} |t|}\int_{-\8}^{t}\|e^{-A^{\lam_1}_s(t-\tau)}p_\p(\tau,\gamma^y_{p,\lam_1}(\tau))-e^{-A^{\lam_2}_s(t-\tau)}p_\p(\tau,\gamma^y_{p,\lam_2}(\tau))\|d\tau\\\\
\leq& k(\p)\int_{-\8}^{t}(t-\tau)^\a e^{-\frac{\eta}{4}(t-\tau)}\(e^{-\frac{\eta}{2} |\tau|}\|\gamma^y_{p,\lam_1}(\tau)-\gamma^y_{p,\lam_2}(\tau)\|_\a\) d\tau\\
&+k(\p)M\int_{-\8}^{t}(t-\tau)^\a e^{-\frac{\eta}{4}(t-\tau)}\big|1 -e^{-(\lam_2-\lam_1)(t-\tau)}\big|d\tau\end{split}\end{equation}
and \begin{equation}\label{eq7}
\begin{split}&e^{-\frac{\eta}{2} |t|}\int_t^{\8}\|e^{-A^{\lam_1}_u(t-\tau)}p(\tau,\gamma^y_{p,\lam_1}(\tau))-e^{-A^{\lam_2}_u(t-\tau)}p(\tau,\gamma^y_{p,\lam_2}(\tau))\|d\tau\\\\
\leq& k(\p)\int_t^{\8}e^{\frac{\eta}{4}(t-\tau)}\(e^{-\frac{\eta}{2} |\tau|}\|\gamma^y_{p,\lam_1}(\tau)-\gamma^y_{p,\lam_2}(\tau)\|_\a\) d\tau\\
&+ k(\p)M\int_t^{\8}e^{\frac{\eta}{4}(t-\tau)}\big|1 -e^{-(\lam_2-\lam_1)(t-\tau)}\big|d\tau.\end{split}\end{equation}
By \eqref{eq5}, \eqref{eq6} and \eqref{eq7}, we derive that \begin{equation}
\begin{split}
&e^{-\frac{\eta}{2} |t|}\|\gamma^y_{p,\lam_1}(t)-\gamma^y_{p,\lam_2}(t)\|_\a\\
\leq&e^{-\frac{\eta}{2}|t|}\int_{0}^t\|e^{-A^{\lam_1}_c(t-\tau)}p(\tau,\gamma^y_{p,\lam_1}(\tau))-e^{-A^{\lam_2}_c(t-\tau)}p(\tau,\gamma^y_{p,\lam_2}(\tau))\|d\tau\\
&+e^{-\frac{\eta}{2} |t|}\int_{-\8}^{t}\|e^{-A^{\lam_1}_s(t-\tau)}p(\tau,\gamma^y_{p,\lam_1}(\tau))-e^{-A^{\lam_2}_s(t-\tau)}p(\tau,\gamma^y_{p,\lam_2}(\tau))\|d\tau\\
&+e^{-\frac{\eta}{2} |t|}\int_{t}^\8\|e^{-A^{\lam_1}_u(t-\tau)}p(\tau,\gamma^y_{p,\lam_1}(\tau))-e^{-A^{\lam_2}_u(t-\tau)}p(\tau,\gamma^y_{p,\lam_2}(\tau))\|d\tau\\
\leq&k(\p)\int_0^{\8}(2+t^{-\a})e^{-\frac{\eta}{4}t}dt\.\sup_{t\in\R}e^{-\frac{\eta}{2} |t|}\|\gamma^y_{p,\lam_1}(t)-\gamma^y_{p,\lam_2}(t)\|_\a\\
&+k(\p)M\int_0^\8(2+t^{-\a})e^{-\frac{\eta}{4}t}\(e^{(\lam_2-\lam_1)t}-1\)dt. \end{split}
\end{equation}
It follows that
\begin{equation*}
\begin{split}\|\xi_p^{\lam_1}(y)-\xi_p^{\lam_2}(y)\|_\a&=\|u_{\lam_1}(0)-u_{\lam_2}(0)\|_\a\\
&\leq\sup_{t\in\R}e^{-\frac{\eta}{2} |t|}\|\gamma^y_{p,\lam_1}(t)-\gamma^y_{p,\lam_2}(t)\|_\a\\
&\leq \frac{k_1(\p)M}{1-M_\p}\int_0^\8(2+t^{-\a})e^{-\frac{\eta}{4}t}\(e^{(\lam_2-\lam_1)t}-1\)dt\\
&\leq \frac{k_1(\p)M}{1-M_\p}\int_0^\8(2+t^{-\a})t \,e^{-[\frac{\eta}{4}-(\lam_2-\lam_1)]t}dt\.|\lam_1-\lam_2|,\end{split}
\end{equation*}where the differential mean value is applied to $e^{(\lam_2-\lam_1)t}-1$ to get the last term.
It is clear that the integral $$\int_0^\8(2+t^{-\a})t \,e^{-[\frac{\eta}{4}-(\lam_2-\lam_1)]t}dt=\int_0^\8(2t+t^{1-\a}) \,e^{-[\frac{\eta}{4}-(\lam_2-\lam_1)]t}dt$$ converges. Therefore $$\xi_p^{\lam_1}(y)-\xi_p^{\lam_2}(y)\|\leq L_2|\lam_1-\lam_2|,$$ where $$L_2:=\frac{k_1(\p)M}{1-M_\p}\int_0^\8(2t+t^{1-\a}) \,e^{-[\frac{\eta}{4}-(\lam_2-\lam_1)]t}dt,$$ and thus $\xi_p^\.(y)$ is Lipschitz continuous on $I_{\lam_0}(\eta/8)$ uniformly on $p\in P$ and $y\in \overline{\mb{B}}_{X^\a_c}(\varrho)$. $\Box$

\end{document}